\documentstyle{amsppt}
\magnification=1200
\nologo

\def\ss{\smallskip}
\def\ms{\medskip}

\def\pa{\partial}
\def\de{\delta_{\pa\O}(x)}
\def\T{\Theta}
\def\L{\Lambda}
\def\oonL{\mathop{{\hbox{$L^{m,2}$}\kern -20pt\raise7pt
\hbox{$\circ$}}}}

\def\O{\Omega}

\def\emp{\emptyset}
\def\bigtimes{\mathop{{\lower2.95pt\hbox{$\wedge$}\kern-6.666pt\raise2.95pt
\hbox{$\vee$}}}\limits}
\def\crd{\cr\noalign{\vskip4pt}}


\def\eq{\eqalign}

\def\e{\epsilon}

\def\l{\lambda}
\def\R{${\text{\bf R}}^N$}
\def\no{\noindent}
\def\bs{\bigskip}
\def\G{\Gamma}
\def\g{\gamma}
\def\p{^\prime}
\def\di{\text{\rm diam\ }}
\def\b{\beta}
\def\a{\alpha}
\def\n{\nabla}

\def\da{d_{\partial \O}(x)}
\def\d{\da}
\def\pn{\par\noindent}

\let\hacek=\v
\def\v{\vert}

\hsize=5in
\vsize=7.6in
\baselineskip=16pt
\def\de{\delta_{\pa\O}(x)}

\def\Wmps0O{W^{m,p}_0(\O,\d^s)}
\def\Wmpt0O{W^{m,p}_0(\O,\d^t)}

\NoBlackBoxes
\NoRunningHeads
\widestnumber\key{ABCDEFG}

\topmatter
\title
Hardy and Hardy PDO type Inequalities in Domains
\\
Part I
\endtitle
\author 
Andreas Wannebo
\endauthor
\abstract
Here is given mainly an extensive treatment of Hardy inequalities in
domains in $\text{\bf R}^N$. 
Part II will continue on this theme and furthermore also treat Hardy PDO 
inequalities.
The latter are Hardy inequalities involving Partial Differential Operators 
instead of gradients in the inequalities. 
General subsets of the function space are also treated. 
Examples are nonnegative cones and are given a special treatment.
\par
A part of the material is given in an encyclopedic fashion. This in
order to get better overview and also be helpful in cases of
outside applications.
\par
These papers are a continuation of [WAN5], which should be regarded a 
prerequisite.
\par
Part II will contain more about Hardy inequalities and will treat
Hardy PDO type inequalities in domains.
\endabstract

\endtopmatter
\document

\centerline{\smc Contents}
\bs

{\bf 0.} Some initial definitions
\ss

{\bf 1.} On the history of this paper
\ss

{\bf 2.} Somewhat on Polynomial Capacities
\ss

{\bf 3.} Hardy and Hardy PDO type inequalities in domains -- some background
\ss

{\bf 4.} On applications
\ss

{\bf 5.} Overview of the results in Part I
\ss

{\bf 6.} The main body of results and proofs -- Hardy inequalities
\bs
\ms

\heading
0. Some initial definitions
\endheading

Let $u$ be a function defined on a subset of $\text{\bf R}^N$
and let $\n^mu$ denote the $m$-th gradient of $u$. This is the set of
partial derivatives in some sense of $u$ of order $m$.
This set written in the multiindex notations is
$\{D^{\a}u\}_{|\a|=m}$.
The $L^p$-norm of this $m$-th gradient can be given in various ways,
usually equivalent. We give the following

$$
||\n^mu||_{L^p}
=(\int|\n^mu|^pdx)^{1\over p}=(\sum_{|\a|=m}\int|D^{\a}u|^pdx)^{1\over p}
\leqno(0.0)
$$
and
$$
\sum_{|\a|=m}(\int|D^{\a}u|^pdx)^{1\over p}.
\leqno(0.1)
$$

Which definition to use is usually only of interest when determining 
e.g. the constants in the inequalities.
The definition is then chosen to be correct with respect to the
application.

The formula (0.1) emphazises the point of view of the norm as a 
sum of seminorms.
\ss

Two versions of Sobolev spaces are the most common.

\proclaim {Definition 0.0}
Let $\O$ be open in \R. 
Let $W^{m,p}(\O)$ be defined as the completion of $C^m(\O)$ in the 
Sobolev norm $||\ ||_{W^{m,p}(\O)}$, which is defined as

$$
||v||_{W^{m,p}(\O)}=\sum_{k=0}^m||\n^kv||_{L^p(\O)}.
$$

If the completion is taken with respect to $C^{\infty}_0(\O)$ instead then
the resulting Sobolev space is denoted $W^{m,p}_0(\O)$.
\endproclaim

These Sobolev space functions are quasicontinuous. This
property means that if the indices for the Sobolev space are $m,p$, 
then such a function is continuous except for an open set of capacity
$\e$ any $\e>0$.
The capacity is the Sobolev space capacity, see [WAN5], here written $C_{m,p}$.
It is equivalent to Bessel capacity $B_{m,p}$  if $p>1$.
The Sobolev functions are defined up to capacity zero.
The capacity $C_{0,p}(M)$ is Lebesgue measure for every $p>0$.
\bs

\heading
1. On the history of this paper
\endheading

The author's research in the area began with a suggestion by L.I. Hedberg
to study the note [ANC1] by Ancona. Ancona studied the following question, when
is a Sobolev space of type $W^{m,p}_0(\O)$ generated by its
nonnegative cone,  i.e. when is every function in the space equal to
the difference of two nonnegative functions in the space? 

The aim was first to improve on the results by Ancona.
However it was found that sometimes his results were the best possible!

This question is connected to Hardy inequalities in domains. 

Hence the study turned into a deeper study of these.
The choice was also motivated by that these inequalities in many cases
are very useful tools and hence progress here 
will influence other areas, such as parts of analysis and mathematical 
physics through their influence on the study of PDE:s, eigenvalues, 
analysis on manifolds, etc. 
\bs

There is also a related question on the compactness and the weights/domains
that make the corresponding imbedding compact. This is not touched upon
here. Certainly it is kind of a twin problem though.

The study of Poincar\'e inequalities for functions in cubes and 
polynomial capacities has been used as a part of this work.
We refer to [WAN5] for this.
The study of these polynomial capacitities has also other applications
to Sobolev space theory. Some aspects of these are treated in [WAN5]
as well.
\ms

In order to illustrate the discussion below on Hardy inequalities for
domains we give a simple formulation of
the Hardy inequality in a domain $\O$ of \R\ namely
$$
\int_{\O}|u|^p\d^{s-mp}dx
\le
A_0\int_{\O}|\n^mu|^p\d^sdx.
\leqno(1.0)
$$
Here $\d$ is the usual distance function from a point $x\in \O$ to
the boundary, $\pa\O$.

The most standard question concerning this inequality is that given
$m,p,s$, which sufficient (necessary) conditions are there on 
$\O$ for inequality (1.0) to hold say for all $u\in W^{m,p}_0(\O)$
with some constant $A_0$.
Furthermore if (1.0) holds then it is also of great interest to get 
information about the best possible constant $A_0$, since it gives 
one of the eigenvalues to a corresponding PDO, e.g. the Laplacian.

This last question is not treated here. The methods
given here are constructive and a value of $A_0$ can be estimated.

-- It should be observed that the particular inequality (1.0) is scale
invariant and this makes it stand out among Hardy inequalities.

In the case of more general weights it should be said that the 
most interesting ones in terms of applications are weights that
are functions of the distance to the boundary. The treatment here 
is much adjusted to this.
\bs

The results given here, say for the simple situation (1.0),
is given in terms of a polynomial capacity. This way both largeness
and shape! of $\O^c$ is measured locally and at all scales.
The bigger this capacity uniformly, the smaller the value of $A_0$.

This gives a formulation of a 
sufficient condition that holds in the case (1.0) if $s<s_0$. Here
$s_0>0$ is calculable.

If $s\ge s_0$ this condition is not enough and a more involved
condition is used. The result is that there is a cost,
a lesser good polynomial capacity has to be used. How much so is regulated
by the badness of  $\pa\O$.
This badness is measured by the value of a dimension of $\pa\O$. The
dimension is designed to reflect local properties uniformly at all scales.

This dimension is denoted $\dim_{loc}$. It was studied by
the author in the mid 80:ies. (There is a typed manuscript by the
author from about that time treating Hardy inequalities.)

The polynomial capacities were originally given by Maz'ya, see [MAZ1],
and also treated in his book [MAZ2]. 
In [WAN5] we gave a more extensive treatment and also constructed a different
polynomial capacity. In [WAN5] also the setup is generalized.

Since the polynomial capacities can estimated from below by
Bessel capacities the same sufficient conditions holds with
Bessel capacities.

Generally this estimate give fewer cases of possible $\O$:s then the
polynomial capacity condition. However these conditions coincide is
when $m=1$ and $p>1$.
\bs

To begin with my technique for these problems only worked for $s<0$,
($m,p$ general).

However Ancona then visited Sweden (Uppsala) and then I 
asked him the question about the case $s=0$ for inequality
(1.0) with $m=1$ and $p>1$ and with this uniform Bessel capacity
condition as discussed above.

Since the question of possible $\O$:s in (0.1) is  in a way is harder 
when $s$ is increasing, this was the first simple but general 
case unknown to me.

Ancona later in [ANC2] answered this question affirmitively and gave a theorem
covering this question just as we discussed. He restricted himself to
$p=2$ and proved the converse in two dimensions.

Then Lewis [LEW] got interested in this part of Ancona's paper. 
He devoted a paper to (1.0) above and 
were able to generalize the statements by Ancona to $m=1$, $p>1$ and $s\le 0$. 
He also proved a converse statement for the cases  with $s=0$ and $p=N$,
(here $m=1$).

The present paper as well as e.g. [WAN2] or [WAN1] 
contains the results of Ancona and Lewis and much more. 
The converse statements by Ancona-Lewis are not treated though.

However the author made a kind of such announcement in [WAN4].
\ms

The methods of the author, Ancona and Lewis are all different, but
the formulation used was given and used first by the
author, then communicated to Ancona and then indirectly to Lewis
through [ANC2].
\bs

There have been several manuscripts by the author with this kind of material
from the mid 80:ies and on. The names of these have been different but
the contents have been similar but expanding.

-- These manuscripts have been circulated.

The version [WAN1], which was not at all the first one,
has been much circulated and dates from (Febr) 1992.
Though also this manuscript has been somewhat expanded later and
circulated as well.

The present paper follows much the manuscript [WAN1] though the part
on polynomial capacities and Poincar\'e inequalities from [WAN1] 
has already been taken up in [WAN5]. 
\bs

The thesis [NYS] by Nystr\"om contains some on these matters.
However his only contribution to the area is 
an estimate of the Maz'ya polynomial capacity 
in the case of domains with boundaries that have the so called Markov 
property. (This will be discussed elsewhere.)

This property is used (by the Ume\aa\ group) as a means to get the same
properties other fractal sets as selfsimilarity gives  this offer.

\heading
2. Somewhat on Polynomial Capacities
\endheading

The first polynomial capacities were invented by Maz'ya. In [WAN5],
(and in [WAN1]) a
deeper treatment is made, which included a generalization and further 
properties. 
Also another and different kind of polynomial capacity was
constructed.

It is the decided opinion of the author that these polynomial
capacities have an important role to play in the theory
related to Sobolev spaces.
The work by the author have so far given a some corroboration.

However we leave this issue to future research to tell.

Since it is rather complicated to go through the matter of definitions,
properties etc. of these two polynomial capacities we refer the reader 
to [WAN5].
This paper is thus a necessary prerequisite for reading the present paper.
(However the reader is refered first to the first part of Section 6
in order to get advice how to get a simpler way to 
a first understanding of these two papers.)

Anyway we now repeat some of the background of polynomial capacities 
in a simplified manner for the benefit of the reader.

These polynomial capacities are based on functions in a unit cube. 
The choice of a cube is merely a practical matter.
Say that someone wants to improve the constants in the Hardy
inequalities that follows from calculations based on the proofs
given here. Then there can be a point to for instance choose 
a unit ball instead.
However in order to have a relationship to the Poincar\'e inequalities,
which is the at issue here, you need that what in [WAN5] 
is called weak Poincar\'e inequalities should hold for the domain.

If one leaves this main
road then of course it is possible to make estimates in some cases. 
But if one treat domains in general and make conditions in terms of
corresponding polynomial capacities then the condition 
may be without meaning.

-- The problem of statements (theorems) and their meaning/contents or
lack of this is a common problem in the area of Sobolev space theory
and requires a constant attention.
\bs

We return to the setting of polynomial capacities in the framework of
a unit cube, $Q$.
We study the simple inequality
$$
||u||_{L^p(Q)}\le C_0||\n^mu||_{L^p(Q)},
\leqno(2.0)
$$
where $u\in {\Cal A}$ and ${\Cal A}$ is a subset of $W^{m,p}(Q)$. 
We assume that the constant $C_0$ is the best possible.
The treatment of general subsets $\Cal A$ given here comes from [WAN1].

Then $C_0^{-p}$ is equivalent to a certain polynomial capacity when $C_0$ is
big enough or as well if the polynomial capacity is small enough.

In this case the polynomial capacity can be chosen as the original polynomial
capacity by Maz'ya. It is here denoted by $\G({\Cal A})$ with suitable indices 
added.

However also the polynomial capacity given by the author also gives a correct
answer. It is denoted by $\T({\Cal A})$. Hence these two types 
of polynomial capacities are here equivalent in the case (2.0)
when $\Cal A$ is varied.	
These polynomial capacities with their respective indices written out 
are in this case (2.0)
$\G_{m,m-1,p}({\Cal A})$ and $\T^{\a}_{m,m-1,p}({\Cal A})$.
Here $m$ is the order of the gradient, $m-1$ is the degree of the
polynomials in the zero space of the RHS and $p$ is exponent in the
norm of the highest order term in the RHS.
(Finally $\a$ is a kind of dummy parameter needed only for 
some proof procedures.)
\bs

Next we study two somewhat more complicated inequalities,

$$
||u||_{L^p(Q)}\le C_1(||\n^{k+1}u||_{L^p(Q)}+||\n^mu||_{L^p(Q)})
\leqno(2.1)
$$

and

$$
||u||_{L^p(Q)}\le A_0||\n^{k+1}u||_{L^p(Q)}+C_2||\n^mu||_{L^p(Q)}.
\leqno(2.2)
$$

In the last case the constant $A_0$ is assumed to be fixed and big enough.

\no
-- Just as before the quantity $C_i^{-p}$ is studied. 

The inequalities hold generally when the respective polynomial
capacity is used and is nonnegative.
As before there is the same relationship between best possible
constant $C_i$ and the respective polynomial capacity.
\ms

In (2.1) the polynomial capacity is denoted $\G_{m,k,p}({\Cal A})$.
\ss

In (2.2) the the polynomial capacity is denoted $\T^{\a}_{m,k,p}({\Cal A})$.
\vskip 0.1truecm

The latter polynomial capacity was (as said before) constructed by the author,
see [WAN5] (or [WAN1]).

In both cases the index $k$ denotes 
that the zero space of RHS is the polynomials of degree less
than or equal to $k$.

A type of inequality of the kind 
(2.2) was first studied by Hedberg [HED]. He gave
a Bessel capacity estimate of the constant which gives a sufficient
condition.With the polynomial capacity $\T$ there is a neccessary and
sufficient formulation instead. It is important to
emphasize the very different nature of the polynomial capacities.
Even in a standard setting. There are key geometrical dependences
that does not exist at all in ordinary (say Bessel) capacities.

-- Anyhow Hedberg managed to get what he wanted, using his estimate.
\bs

If a situation is hand that makes both
corresponding polynomial capacities equivalent, 
then of course the inequality (2.2) is better than (2.1). 

-- However there is reason! to pose a question,
which is not natural because it seems contraintuitive.

Anyway experience give some hope for an affirmative answer anyway.

\proclaim{Open question 2.1} 
\ss

\no
-- Are these two types of polynomials capacities equivalent?

\no
-- Refrased, does (2.1) imply (2.2)? Maybe with new fixed constants.
\endproclaim
\bs

\no
{\bf Comment on vocabulary.}
\ss

From our point of view the word polynomial in the expression 
polynomial cpapacity here
should be regarded as derived from 
that the set of {\bf polynomials} of degree
less or equal to $m-1$ 
is the zero space of the seminorm $||\n^mu||_{L^p(Q)}$.
\ss

The use of the word {\bf capacity} should {\bf not}
be seen as derived from the word capacity as used in any ordinary sense.
Instead it should be seen as the fact that the polynomial capacities 
are intimately related to say Sobolev (or Bessel potential) 
space capacities in force of the actual formulas.
\ss

-- Formally, i.e. in a logical sense, the concept of polynomial capacity
is entirely different from that of capacity.

\heading 
3. Hardy and Hardy PDO type Inequalities in Domains -- Background
\endheading 

The first person to generalize Hardy inequalities to more general domains
was Ne\hacek cas. 
He took the standard case in one dimension
to the corresponding case for Lipschitz domains. 

Then Kufner used Ne\hacek cas' ideas and generalized to H\"older domains. 
He wrote a book [KUF] with this as the main message, together with
treatment of applications etc.

It is a thin volume but has anyway been influential. In fact since his
results were not contested Kufner thought that maybe his parameter settings
might be optimal, until he got to know our work in the field.

The results in [KUF] on Hardy inequalities in domains and imbeddings
are special cases of the single theorem in [WAN3], part of author's 
thesis 1991. However Kufner gives also another setting (variant) 
of inequalities that we do not discuss.

Again we return to the following aspect.

--In my opinion there is a big problem in this area of Hardy inequalities
and also other areas. Namely sometimes there are results given that
are posing as very good ones -- for instances necessary and
sufficient, but instead these results are dubious since they do 
not lead very far when you want a specific answer.

This can be said about the result in [MAZ2] p.113, i.e. if it is regarded as
a Hardy inequality. He gives a very
general statement about first order Hardy inequalities with general
domains and general measures as weights in the LHS and an expression
for the RHS that can be very involved and can be used for a defining
a capacity of Choquet type.

The problem is that he makes a comparison between the values of
the measure and the capacity every compact in the domain
with respect to the full domain.
The infimum you get is the best constant in the inequality.

The problem is that these numbers usually cannotbe calculated.

Maz'ya discusses the problem in [NIK] p.153 and that shows that he is 
aware the problem. 

Anyway his result in the context of the book [MAZ2] is a very important lemma,
that gives sharp constants in unweighted inequalities like the Sobolev
inequality and variants, see also Stredulinsky [STR] for a different
short proof. It is based on another technique -- introduced by Maz'ya
-- ``capacitary integrals''.

Then we have the many results in [GUR-OPI]. They do not give examples 
and it seems hard/impossible? to get interesting ones. 
This makes gives their results a dubious ring.

Horiuchi in [HOR] gives a useful condition for a H\"older
boundary such that the boundary has dimension less than the 
dimension of the domain minus one.
Here it ought to be worthwhile to generalize along the lines of [WAN3].  

\heading
3.1 Hardy PDO type Inequalities
\endheading

We make a somewhat brief and uncomplete history.
\bs

These inequalities appeared in [WAN1] and shortly thereafter also in
mathematical journals. 

The general idea here have been to use as condition that if we take
a point in $\O$ with a distance $r$ from the boundary then there
should exist a ball with radius $\l r$ for some fixed $\l$ such that
it within this ball exist a ball in the complement of $\O$ with 
radius $\mu r$ with $\mu$ fixed. Then the situation should about the
same as for the corresponding Hardy inequality.

In fact in [WAN1] there has much more results than those that appeared
in journals. However this part of [WAN1] needs some reorganization.
\bs

Furthermore the papers [SHA] and [H-K-S] contain results which
seem to be very much related the question of Hardy PDO type
inequalities, though the connection has not been worked out yet.
They have many applications of their results.
\ss

A treatment of Hardy PDO type inequalities will be included in Part II. 

\heading
4. On Applications
\endheading

The Hardy inequalities for domains and the accompanioning 
problem of compactness of
the corresponding inequality (which is not treated here) are heavily
linked to many areas of mathematical physics and more generally
the study of PDE:s as well as their eigenvalue problems. However the
potential use is even greater.

-- Another aspect is that these Hardy inequalities can be seen as models
for other kinds of inequalities.

One main aspect is that the Dirichlet problem generalizing Poission's 
equation can be treated for very bad domains using a variation of
Lax-Milgram's lemma. For this see the relevant chapter in [KUF].
In fact the results there are extremely general.

The theory given there can be said to be just waiting for better 
Hardy inequalities in order to be upgraded,
also the book [KUF-S\"AN] by Kufner and S\"andig can be recommended as source
of many examples of applications.

-- Ancona's paper [ANC2] is quite helpful and with an alternative point of
view to that in the Kufner book [KUF]. He treats the classical
Dirichlet problem more classically.

Nystr\"om in [NYS] discusses some aspects of the Dirichlet problem.

Horiochi have motivated his study of Hardy inequalities for domains by
the fact that they are suited for treating Dirichlet problems with
pertubed ellipticity.
This aspect is also treated in [KUF-S\"AN].

A field in physics where these Hardy inequalities are much
present is General Relativity. 

Generally when PDE:s are treated in physics there is often a possible use 
for Hardy inequalities.

Also generally speaking when there is progress on Hardy inequalities in
domains in \R, then this can be translated to manifolds with/(without) 
boundary.

A special question here, where Hardy inequalities in domains
are of interest is the so called Yamabe problem, see Maz'ya's opinion
on this, [NIK].

\heading
5. Overview of the Results
\endheading

{\bf I.}
\ss

The first result of the paper is a summation lemma, Lemma 6.4. 
It is used to go from local information in the sense of local
inequalities to a global one. 
It involves sums of integrals over
enlarged Whitney cubes and concerns general functions. 
Since it does not involve
Sobolev functions or special kind of functions it can be used in many
contexts.
\bs

{\bf II.}
\ss

 We define a dimension here denoted $\dim_{loc}$, see Definition 6.8. 
It appears naturally by the use of local H\"older inequalities.

We also define another dimension $\dim_{mc,loc}$ using a similar
definition which is based on Minkowski content, see Definition 6.9. 
The dimensions can be proved to be equal for compacts but the
proof is not included.
Also the same can be asked about a similar constructed dimension
based on Hausdorff measure instead.
\bs

{\bf III.}
\ss

We make a minor variation of the definition of H\"older quotient
in order to suit the summation process, see Definition 6.12.
\bs

{\bf IV.}
\ss

The main result on Hardy inequalities is collected in Theorem 6.17.
This theorem should be regarded as a kind of look-up table, since there
are too many possibilities squeezed into one theorem to make any nice
version.

The results are given with the sum of two integrals in the LHS since
this sometimes is advantageous. Such situations will be treated in Part II.

The two different types of polynomial capacities 
make a difference when the two-integral formulation is used.

The reader is adviced to look through the passage ``Advice for the
first reading'' given in the beginning of section 6 before more 
careful reading.

The contents of Theorem 6.17 follows in more detail from the knowledge of 
the respective capacities, i.e. see [WAN5]. 
\bs

{\bf V.}
\ss

Corollary 6.19 is a rewriting of the main theorem in certain
situations for the benefit of the reader. Here also [WAN5] is used.
We observe for instance that if $m=2$ then nonnegative functions often have 
better 
$\T$- capacities and that way they get better Hardy inequalities as well.

Corollary 6.19 
contains a long list but is a god place to look for adequate information.
\bs

{\bf VI.}
\ss

Theorem 6.20 gives sufficient conditions for when certain
Sobolev spaces equal the difference formed from its nonnegative
cone. 
Let lower index $+$ denote the non-negative cone. Then this is
written

$$
W^{m,p}_0(\O,\d^s)=W^{m,p}_0(\O,\d^s)_+-W^{m,p}_0(\O,\d^s)_+,
\leqno(5.0)
$$

i.e. also weighted Sobolev spaces are covered.
Actually the  conditions given make substantial improvement to the
corresponding theorem by Ancona in [ANC1]. He studies the cases 
with $s=0$, $p>1$ and $\O$ Lipschitz, or $p>N$. 

Theorem 6.23 is just an  upgrading of his result with essentially his
original proof.

Theorem 6.21 treats a similar problem and depends on
Theorem 6.23 together with our Hardy inequality results. 

With one of several given conditions given (listed) we have that 
$u\in W^{2,p}_0(\O,\d^s)$ implies 

$$
u\in W^{2,p}_0(\O,\d^s)_+-W^{2,p}_0(\O,\d^s)_+
\leqno(5.1)
$$

if and only if 

$$
\int|u|^p\d^{-2p+s}dx<\infty.
\leqno(5.2)
$$
\bs

{\bf VII.}
\ss

We formulate a rather spectacular conjecture.
\bs

\proclaim{Conjecture 5.0} 
\ss

Let $m$ be odd, $p>1$ then 

$$
W^{m,p}_0(\O)=W^{m,p}_0(\O)_+-W^{m,p}_0(\O)_+
\leqno(5.3)
$$

holds for all $\O$ open in \R.
\ss

Let $m$ be even and positive then there always exist some $N,p,\O$ ,
such that (5.3) does not hold.
\endproclaim

\proclaim{Remark 5.1} 
\ss

It be conjecured that the weighted version of Conjecture 5.0 would hold.

This is not done since then features are introduced that makes the answer
possibly out of reach. 
\par
-- For ever!
\endproclaim

\heading
6. The main body of results and proofs
\endheading

\proclaim{Advice for a first reading}
\endproclaim

The nature of the information given here and in [WAN5] certainly
motivates some advice.  The actual formulas for the main results are
given as inequalities in formulas (6.24) and (6.25). Hence LHS and RHS 
below refer to these.
\ss

(i) The first term in the RHS, which involves $p_1$,
usually is used only in special situations.
If this term is skipped it simplies.
To begin with there is no longer any difference then between the $\G$- and
the $\T$-capacities and this reduces the number of cases.
Furthermore Lemma 6.14 can be postponed.
\ss

(ii) The key lemma is Lemma 6.4. -- It should be read.
\pn
\ \ \ \ \ \ \ \ \ The main theorem is Theorem 6.17. 
\ss

(iii) Temporarily the information on polynomial capacities given in
[WAN5] can wait.
-- Instead they can for the moment be seen as a representation of the 
constant in the Poincar\'e inequality in a cube, see Section 2.
\ss

(iv) Put $q=p$ and get some simplification.
\ss

(v) Now only Case A -- (6.25) and Case E are left. 

Some part of the proof is given in Case A -- (6.24).
Hence it is needed to read this too. 
Choose one of the two suggested cases.
\ss

(vi) To understand the outcome. Some understanding of the polynomial 
capacities is needed.
To faciliate go to Corollary 6.19 first. 
Then choose a special case of interest
and in accordance to the reduction made.
Then return to [WAN5], look for the fact and its
explanation.

\proclaim{Observe} When constants $A$ are used several times their
values may change as the steps in a proof proceed.
\endproclaim

To begin with we first introduce the Whitney cube concept.
It is in fact a way to look upon the open sets (with non-empty complement)
as combinatorial objects.
\ss

\proclaim {Definition 6.0} 
\ss

Given an open set $\O\in\text{\bf R}^N$ and
a set of cubes ${\Cal F}_{\O}$. This set is called the Whitney cubes 
of $\O$ or are said to form a Whitney decomposition of $\O$, if
the following conditions hold.

The cubes in ${\Cal F}_{\O}$ are open, dyadic and disjoint.
Furthermore
$$
\O=\bigcup_{Q\in{\Cal F}_{\O}}\bar Q,
\leqno(6.0)
$$
$$
\di Q\le dist(\O^c, Q)\le 4\di Q
\leqno(6.1)
$$
and, if $\bar Q\cap\bar Q'\not=\emp$, then
$$
{1\over 4}\le {\di Q\over \di Q'}\le 4.
\leqno(6.2)
$$
\endproclaim

\proclaim{Theorem 6.1} (Whitney.)
\ss

A set of Whitney cubes as defined above exist (non-uniquely) for every
open subset of \R with a non-empty complement.
\endproclaim

Proof. See e.g. the book by Stein, [STE].
\ss

We now give a very useful Summation Lemma involving integrals and weights.

-- This lemma gives a generic method for many situations.

First some notation.

\proclaim{Definition 6.2}
\ss

Let $\O$ be open \R. Let ${\Cal F}_{\O}$ be a Whitney decomposition of $\O$.
For $Q\in {\Cal F}_{\O}$ define a new cube as follows. Take any
point $x_0$ on $\pa\O$ that has smallest distance to $Q$. Then let
$x_0$ be centre of a cube $R_Q$ which has the smallest side length and
covers $Q$.
\endproclaim
\ms

\proclaim{Notation 6.3}
\ss
Fix $Q\in{\Cal F}_{\O}$.
Denote with $\cdot\tilde{\phantom{o}}$ for the scaling which has the
property that $R_Q$ is mapped on a unit cube. This cube is denoted
$\tilde R_Q$. By the same scaling $\O$ 
is mapped on $\tilde{\O}$
etc. When $\cdot\tilde{\phantom{o}}$ is
used it should be clear what $Q$ is refered to.
\endproclaim

\proclaim{Lemma 6.4} 
\ss

Let $\O$ be an open set in \R\ with
a Whitney decomposition ${\Cal F}_{\O}$.
Let $f$ be a nonnegative function on \R\ with $f|_{\O^c}=0$ and 
let $s>0$. Then
$$
\sum_{Q\in{\Cal F}_{\O}}(\di Q)^{-s}\int_{R_Q}f\d^sdx
\le {A(N)\over 1-2^{-s}}\int_{\O}fdx.
\leqno(6.3)
$$
\endproclaim

Proof.
\ss

It follows from the properties of Whitney cubes that 
$d(x,\pa\O)\le 5\di Q$ for $x\in Q$.
This together with a change of order of summation gives that

$$
\eq{
&\sum_{Q\in{\Cal F}}(\di Q)^{-s}\int_{R_Q}f\d^sdx
\crd
&\le
\sum_{Q\in{\Cal F}}\sum_{Q\p\in{\Cal F}\atop Q\p\cap R_Q\not=\emp}
(\di Q)^{-s}(5\di Q\p)^s\int_{Q\p}fdx
\crd
&=\sum_{Q\p\in{\Cal F}}\sum_{Q\in{\Cal F}\atop Q\p\cap R_Q\not=\emp}
5^s(\di Q)^{-s}(\di Q\p)^s\int_{Q\p}fdx.
\cr}
\leqno(6.4)
$$

We want to evaluate the inner sum of the RHS of (6.4) for fixed $Q'$.
We have to estimate some of the entities involved.
First we estimate $\di Q'$, when $Q'$ intersects $R_Q$.
\ms

To this end we use a sublemma.
\ms

\proclaim{Sublemma 6.5}
\ss

With the situation as in Lemma 6.4, 
there is a constant $a$ such that 
$$
{\di Q'\over\di Q}>a
\leqno(6.5)
$$
implies that $Q'\cap R_Q$ is empty. The value of $a$ can be taken as 
$a=5\sqrt N$. 
\endproclaim
\ss

Proof of Sublemma 6.5.
\ss

First we determine how large $\di R_Q$ can be.

The quotient $\di R_Q/\di Q$ takes its largest value if $Q$
is situated exactly at the middle of a face of $R_Q$ and 
$4\di Q=\text{dist}(Q,\pa\O)$,
since if we move $Q$ 
about in $R_Q$ following a face we get an equal or longer distance
otherwise and 
then $Q$ has to be made larger in order to satisfy the Whitney cube condition
$4\di Q\ge \text{dist}(Q,\pa\O)$. 

Hence, the side of $R_Q$ has length at most $10\di Q$ and the
diameter of $R_Q$ has length at most $10\sqrt N\di Q$.
Now according to the Whitney cube property $\di Q\p\le 5\sqrt{N}\di Q$.
\ss

End of proof of sublemma 6.5.
\ms

The proof of Lemma 6.4 continues.
\ss

The diameters of the Whitney cubes are dyadic and
there is no restriction to take $\di Q=2^{-n}$ and $\di Q\p=2^{-k}$, with
$k$ and $n$ integers. By 
Sublemma 6.5 we have that 
$k+{}^2\log (5\sqrt N)\ge n$ implies that $Q'\cap R_Q=\emp$.

It follows from a volume consideration that
there are at most only a fixed number $A(N)$ of
$Q$:s with $n$ fixed, such that $R_Q$ intersects $Q'$.
If $Q'\cap R_Q\not=\emp$, then the cube $Q$ lies in a ball with 
centre in the centre of $Q'$ and 
$$
\text{radius of ball}=\di R_Q+(1/2)\di Q'=10\sqrt N \di Q+(5/2)\sqrt N \di Q.
$$
The conclusion follows.
\ms

Now we can evaluate the inner sum in (6.4)

$$
\eq{
\sum_{Q\in{\Cal F}\atop Q\p\cap R_Q\not=\emp}
5^s(\di Q)^{-s}(\di Q\p)^s
&\le
A(N)\sup_k\ \{\sum_{n=-\infty}^{k+\lfloor{}^2\log (5\sqrt N)\rfloor}
2^{(n-k)s}\}
\crd
&\le 
{A(N)\over 1-2^{-s}}.
\cr}
\leqno(6.6)
$$
\ss
Lemma 6.4 now follows from (6.5) and (6.6).
\ss

End of proof of Lemma 6.4.
\ms

{\bf Observe:} 
\ss

It holds that for $s$ small

$$
{1\over 1-2^{-s}}\sim {1\over s}.
\leqno(6.7)
$$

\proclaim{Definition 6.6}
\ss

Let $\O$ be open in \R\ and $s$ be real, then define
$$
G_s(\O)=\sup_{Q\in{\Cal F}}(\di Q)^{-N+s}\int_{R_Q\cap\O}\d^{-s}dx,
\leqno(6.8)
$$
where ${\Cal F}_{\O}$ is a Whitney decomposition of $\O$.
\endproclaim

\proclaim{Observation 6.7}
\ss

$G_s(\O)$ is invariant when $\O$ and ${\Cal F}_{\O}$ are dilated.
\endproclaim

We will define a useful related dimension concept.

\proclaim{Definition 6.8}
\ss

Let 
$$
\dim_{loc}(\pa\O,\O)=N-s_0,
\leqno(6.9)
$$
where 
$$
s_0=\sup\{s:G_s(\O)<\infty \}.
\leqno(6.10)
$$
\endproclaim

This concept goes back to our work on the problem of Hardy
inequalities in the mid 80-ties. (Documented in a typed
manuscript from that time.)
\ms

We discuss this dimension concept only as a remark.
The definitions of Hausdorff measure, $h_d$,
and Hausdorff dimension, $\dim_h$, as well the
definitions of Minkowski content $mc_d$ and Minkowski dimension
$\dim_{mc}$ are assumed to be known. 

\proclaim{Definition 6.9}
\ss

Define 
$$
\dim_{mc,loc}(\pa\O,\O)=\inf_d \{d:MC_d<\infty\}
\leqno(6.11)
$$
with
$$
MC_d=\sup_{Q\in{\Cal F}_{\O}}(mc_d(\tilde R_Q\cap\tilde{\pa\O)})
\leqno(6.12)
$$
and define in the same way
$$
\dim_{h,loc}(\pa\O,\O)
\leqno(6.13)
$$
\endproclaim

\proclaim{Theorem 6.10}
\ss

It holds for closed sets that

$$
\dim_{loc}=\dim_{mc,loc}
\leqno(6.14)
$$ 
\endproclaim

Proof. Only in the typed notes from the 80:ies.

\proclaim{Question 6.11}
\ss

Does it hold that

$$
\dim_{loc}=\dim_{h,loc}.
\leqno(6.15)
$$ 
\endproclaim

For some cases in the main theorem we want a variant of H\"older
quotients that performs well under summation.

\proclaim {Definition 6.12}
\ss

Let $\O$ be open in \R, $0<\l \le 1$ and $h$ be nonnegative integer. 
The pointwise H\"older quotient is defined by
$$
||\n^hu||_{H^{\l,pnt}(\O)}=
\sup_{x\in \O}
\ \sup_{|\a|=h}
\ \limsup_{y\rightarrow x}
\{{|D^{\a}u(x)-D^{\a}u(y)|\over |x-y|^{\l}}\}.
\leqno(6.16)
$$
\ss

In weighted formulation with weight function $\kappa$

$$
||\n^hu||_{H^{\l,pnt}(\O,\kappa)}=
\sup_{x\in \O}
\ \sup_{|\a|=h}
\ \limsup_{y\rightarrow x}
\{{|D^{\a}u(x)-D^{\a}u(y)|\over |x-y|^{\l}}\cdot\kappa(x)\}.
\leqno(6.17)
$$
\endproclaim

The following definition is made for examples later.

-- It is not generally agreed on one definition of selfsimilarity.
Here is a definition suitable for the discussion here.

\proclaim{Definition 6.13}
\ss

Let $K$ be a closed set in \R\ and let ${\Cal F}_{\O}$ 
be a Whitney decomposition of 
$K^c$. 

Then $K$ is {\it selfsimilar} if for every $Q\in {\Cal F}_{\O}$, 
$\tilde R_Q$ contains a ball $\tilde B_Q$
such that $\di \tilde B_Q$ is uniformly bounded off zero and every
$\tilde B_Q\cap \tilde K$ has the property that there is a 
similarity transformation taking one to the other. 
\endproclaim

The following lemma is a consequence of a theorem by Sobolev
on equivalent norms in Sobolev space.

\proclaim{Lemma 6.14} 
\ss

Given conditions
\ss

(i) \ Let $Q_0$ be unit cube, $Q$ cube with $Q\subset Q_0$ and 
sidelength $l(Q)=a$,
\ss

(ii) \ let $p\ge 1$ and let $m,k$ with $m>k+1$ be positive integers,
\ss 

(iii) \ let $p_1$ with $0<p_1\le{N p\over N-(m-k-1)p}$ for $N>(m-k-1)p$, 
\ss

(iv) \ let $0<p_1<\infty$ for $N=(m-k-1)p$,
\ss

(v) \ let $0<p_1\le\infty$ for $N<(m-k-1)p$,
\ss

(vi) \ let $\n^mu$ etc. be the vectors of all weak derivatives of
order $m$ of $u$ that are functions a.e.
\ms

Then it holds that there exist a constant $A=A(N,m,p,p_1,a)$
independent of the position of $Q$ in $Q_0$  such that 
$$
(\int_{Q_0}|\n^{k+1}u|^{p_1}dx
)^{1\over p_1}
\le 
A((
\int_Q|\n^{k+1}u|^{p_1}dx
)^{1\over p_1}
+
(\int_{Q_0}|\n^mu|^pdx)^{1\over p}).
\leqno(6.18)
$$
\endproclaim

Proof. It is enough to prove 
$$
(\int_{Q_0}| D^{\a}u|^{p_1}dx)^{1\over p_1}
\le 
A(
(\int_Q| D^{\a}u|^{p_1}dx)^{1\over p_1}
+(\int_{Q_0}|\n^{m-k-1}D^{\a}u|^pdx)^{1\over p})
\leqno(6.19)
$$
for all $\a$, $|\a|=k+1$.
\ss

Hence define
$$
F(D^{\a}u)
=
\min_{Q\subset Q_0\atop l(Q)\ge a}
(\int_Q|D^{\a}u|^{p_1}dx)^{1\over p_1}.
\leqno(6.20)
$$

The functional $F$ is obviously continuous in $W^{m-k-1,p}(Q_0)$.
Furthermore if $P\in {\Cal P}_{m-k-2}$ and $P\not=0$, 
then $F(P)\not=0$. 
These conditions gives that the result follows from 1.1.15 in [MAZ2]. 
\ss

End of proof.
\bs

\centerline{\smc The Main Results on Hardy inequalities}
\bs
\ss

The main Theorem and other results are in order to avoid endless 
repetetivity given as
\bs

\hskip 3.5truecm
{\bf
Statements structured as
}
\vskip 0.3truecm
\hskip 3.5truecm
{\bf -- A Look-up Table With Inputs.
}
\ms

Thus minimizing the formulation of the theorem to about two pages only!
\ms

Later on some of the information in [WAN5] on polynomial and 
usual capacities are used to make lists that are less general
but easier to grasp and overview.
\ms

{\bf Discussion 6.15.}
\ss

The results from now on can be said to be organizied in a somewhat
peculiar manner. 
Furthermore not all situations that can covered by the methods here are
treated. -- A selection has been made.
This selection is done to show different possibilities
when using different extra ideas and also to show their
effects, i.e. the outcome for the possible Hardy inequalities. 

-- By a consistent use of power type weight with respect to distance to the
boundary (in the RHS) a concentration to the in applications most 
important situations is achieved.
-- Certainly the readibility suffers anyway and certainly 
this anyway is in nature of the subject itself.
\ms

In order to simplify the exposition we give some special notation.

\proclaim{Notation 6.16} (For Theorem 6.17.)
\ss

Let $m,k$ be integers and $p\ge 1$ and $p_1>0$.

Let $\O$ be open proper subset in \R, with Whitney decomposition 
${\Cal F}_{\O}$.

Let ${\Cal A}$ be a subset of $W^{m,p}_{loc}$ with ${\Cal A}|_{\O^c}=0$ q.e.

Let $\de$ be the regularized distance function, $\de\sim\d$ and 
$\de\in C^{\infty}(\O)$,
\linebreak 
see [STE].

As shorthand denote
$$
\sum_Q=\sum_{Q\in {\Cal F}_{\O}}.
$$

Denote $[p,p_1]=\max\ \{p,p_1\}$.

Define for $x\in Q\in {\Cal F}_{\O}$

$$
\eq{
\G_{m,k+1,p,{\Cal A}}(x)
&=
\G_{m,k+1,p}({\Cal A}|_{\tilde R_Q\cap\tilde{\O}^c}),
\crd
\T^{\a}_{m,k+1,p,{\Cal A}}(x)
&=
\T^{\a}_{m,k+1,p}({\Cal A}|_{\tilde R_Q\cap\tilde{\O}^c}).
\cr}
\leqno(6.21)
$$

Let

$$
f_{p,p_1}:{\Cal F}_{\O}\rightarrow [0,1].
$$

Let $s$ be given by 
$$
s'=s/(s-1)\text{  and  }s'={\max\ \{p,p_1\}\over q}
$$
and then let $f_{p,p_1}\in l^s$ with norm 1.

Denote $f(x)=f_{p,p_1}(Q)$.
\endproclaim
\ss

Theorem 6.17 is organized as follows. First some general conditions are
given that do not depend on which of the inequalities (6.24) or (6.25)
is at hand nor the ``Cases''. Then two different 
the two types of Hardy inequalities (6.24) and (6.25) are given 
and since they are different there are also 
given general conditions to each of them.
After this the different ``Cases'' are given. 
They include more conditions needed
as well as ``inputs'' to the inequalities (6.24) and (6.25).
\ss

This system of presentation is also used later on.

\proclaim{Theorem 6.17} 
\ss

The results are given in notation 6.16.
\ss

Let $u\in {\Cal A}$.
\ss

In both (6.23) and (6.24) let 

$$
s_1=-(m-k-1)p_1-N+{p_1\over p}(s+N).
\leqno(6.22)
$$

Preconditions for (6.23).
\ms

(i) \ $(m-h)p>N>(m-h-1)p$,
\ss

(ii) \ $0<\l\le m-h-{N\over p}$,
\ss

(iii) \ $0<\l<1$,
\ss

(iv) \ $t=m-h-\l-{s+N\over p}$.
\ms

Preconditions for (6.24).
\ms

(i) \ $0<q\le{pN\over N-mp}$,
\ss

(ii) \ $N>mp$,
\ss

(iii)  \ $0<q<\infty$, when $N=mp$ 

\ \ \ \ \ \ and $q\le\infty$ when $N<mp$,
\ss

(iv) \ $t=mq-({q\over p}-1)N-{q\over p}s$.

$$
\eq{
||\n^hu||_{H^{\l,pnt}(\O,\L(x)^{1\over p}\d^{-t})}
\le&
A((\int_{\O}|\n^{k+1}u|^{p_1}\Lambda_1(x)^{{p_1\over p}}\d^{s_1}dx)
^{1\over p_1}
\crd
&+(\int_{\O}|\n^mu|^p\d^sdx)^{1\over p}).
\cr}
\leqno (6.23)
$$

$$
\eqalign{(\int_{\O}{| u|^q\Lambda(x)^{q\over p}\over\d^t}dx)^
{1\over q}\le&
A((\int_{\O}|\n^{k+1}u|^{p_1}\Lambda_1(x)^{p_1}
\d^{s_1}dx)^{1\over p_1}\crd
&+(\int_{\O}|\n^mu|^p\d^sdx)^{1\over p}).
\cr}
\leqno (6.24)
$$
\ss

Then (6.23) and (6.24) holds respectively when the general condition 
(6.22), the preconditions as well as the conditions in 
Case A-E are satisfied in the form the respective inputs give.
\bs

{\smc Case A.}
\ss

(i) \ Let ${\Cal A}$ be arbitrary,
\ss

(ii) \ let $s<0$, $p\ge 1$ and $\L_1(x)=1$,
\ss

(iii) \ in $(6.23)$ let $\L(x)=\G_{m,k,p}(x)$,
\ss

(iv) \ in $(6.24)$ let $\L(x)=\G_{m,k,p}(x)$ 
      if $q\ge \max\lbrace p,p_1\rbrace$,

\ \ \ \ \ \ and if  $0<q<\max\lbrace p,p_1\rbrace$ then 
$\L(x)=\G_{m,k,p}(x)f(x)^{p\over q}$. 
\ss

Here $A=A(N,m,p,p_1)(1-2^{-s})^{-{1\over p}}\sim A(N,m,p,p_1)s^{-{1\over p}}$
for small $s$.
\bs

{\smc Case B.}
\ss
(i) \ Let ${\Cal A}$ be arbitrary,
\ss
(ii) \ let $\L_1(x)=1$, 
\ss
(iii) \ let $1\le p_0<p$,
\ss
(iii) \ let $\dim_{loc}(\partial\O,\O)<N$,
\ss
(iv) \ let $s<{p-p_0\over p_0}(N-\dim_{loc} (\partial \O,\O))$,
\ss
(v) \ in $(6.23)$ let $\L(x)=\G_{m,k,p_0}(x)$,
\ss
(vi) \ in $(6.24)$ let $\L(x)=\G_{m,k,p_0}(x)$, if 
$q\ge\max\lbrace p,p_1\rbrace$ and $q\le p_0^*$ with $p_0^*$ 

\ \ \ \ \ (Sobolev exponent)

\ \ \ \ \ but if $0<q<\max\{p,p_1\}$, then 
$\L(x)=\G_{m,k,p_0}(x)f(x)^{p\over q}$. 
\ss

Here $A=A(N,m,p,p_1, s, \O)$
\ss

{\bf Comment.} The upper bound on $q$ in (vi) is not made optimal.
\bs

{\smc Case C.}
\ss
(i)\ Let ${\Cal A}$ be arbitrary,
\ss
(ii)\ let $s<0$, 
\ss
(iii)\ let $p\ge 1$,
\ss
(iv)\ let $\L_1(x)=\T^{\a}_{m,k,p}(x)$,
\ss
(v)\ in $(6.23)$ let $\L(x)=\T^{\a}_{m,k,p}(x)$,
\ss
(vi)\ in $(6.24)$ let $\L(x)$ be the same if $q\ge \max\lbrace p,p_1\rbrace$

\ \ \ \ \ \ and if $0<q<\max\{p,p_1\}$, then let 
$\L(x)=\T^{\a}_{m,k,p}(x)f(x)^{p\over q}$.
\ss 

Here $A=A(N,m,p,p_1)(1-2^{-s})^{-{1\over p}}$.
\bs

{\smc Case D}
\ss
(i) \ Let ${\Cal A}$ be arbitrary,
\ss
(ii) \ let $1\le p_0< p$ and let $\text{dim}_{loc}(\partial\O,\O)< N$,
\ss
(iii) \ let $s<{p-p_0\over p_0}(N-\dim_{loc} (\partial \O,\O))$,
\ss
(iv) \ let $\L_1(x)=\T^{\a}_{m,k,p_0}(x)$,
\ss
(v) \ in $(6.23)$ let $\L(x)=\T^{\a}_{m,k,p_0}(x)$,
\ss
(vi) \ in $(6.24)$ let $\L(x)$ be the same if $q\ge \max\lbrace p,p_1\rbrace$

\ \ \ \ \ \ and if $0<q<\max\{p,p_1\}$, then let 
$\L(x)=\T^{\a}_{m,k,p_0}(x)f(x)^{p\over q}$. 
\ss

Here $A=A(N,m,p,p_1,s,\O)$.
\bs
 
{\smc Case E -- $(6.24)$.}
\ss
(i) \ Let ${\Cal A}\subset C^{\infty}_0(\O)$,
\ss

(ii) \ let $\G_{m,m-1,p,{\Cal A}}(x)\ge b>0$ for all $Q\in {\Cal F}$,

\ \ \ \ \ (where the index $m-1$ gives the one-integral case in the RHS),
\ss

(iii)\ let $s_0>0$ be a constant that can be calculated, let $s<s_0$,

\ \ \ \ \ \ let $\L(x)=1$ and $q\ge p$.
\ss

Here $A=A(N,m,p,b)$ and $s_0=s_0(N,m,p,b)$.
\ms

{\bf Comment on Case E.} 
\ss

The situations in the previous Cases with with $p_1$-term, $\T$-,
$\G$-capacities, $\dim_{loc}$, $f$-case, (6.23) and (6.24)
have already been treated and can be adjusted to the Case E as well.
Hence one situation only is treated, also with the simplification
that $\Cal A$ is subset of $C^{\infty}_0(\O)$.
\endproclaim
\ss

\smc{PROOFS:}
\ms

Proof of Case A -- (6.23)
\ss

Obviously

$$
||\n^hu||_{H^{\l,pnt}(\tilde Q)}
\le
||\n^hu||_{H^{\l,pnt}(\tilde R_Q)}
\le
||\n^hu||_{H^{\l}(\tilde R_Q)}.
\leqno(6.25)
$$

By a Poincar\'e inequality for a cube with H\"older seminorms 
in the LHS, see [WAN5],
and (6.25)

$$
\eq{
||\n^hu&||_{H^{\l,pnt}(\tilde Q)}\le
\crd
&\le
{A\over\G_{m,k,p,{\Cal A}}(x)^{1\over p}}
((\int_{\tilde R_Q}|\n^{k+1}u|^{p_1}dx)^{1\over p_1}
+(\int_{\tilde R_Q}|\n^mu|^pdx)^{1\over p}))
\crd
&\le
{A\over\G_{m,k,p,{\Cal A}}(x)^{1\over p}}
((\int_{\tilde Q}|\n^{k+1}u|^{p_1}dx)^{1\over p_1}
+(\int_{\tilde R_Q}|\n^mu|^pdx)^{1\over p})).
\cr}
\leqno(6.26)
$$

Dilate back. The dilation is a change of the independent variable with a
factor $\gamma={\di R_Q\over\di\tilde {R}_Q}$ involved. 
In this construction $\di\tilde{R}_Q=1$.

$$
\eq{
\text{\bf Dilation table}&\text{\ --\ dilation factor $\g$}
\crd
\n^k\tilde u &\rightarrow \gamma^k\cdot\n^ku
\crd
d_{\partial\tilde\O}(x) &\rightarrow \gamma^{-1}\cdot\da
\crd
dx &\rightarrow \gamma^{-N}\cdot dx
\cr}
\leqno(6.27)
$$

Next multiply with $\G_{m,k,p,{\Cal A}}(x)^{1\over p}$.
Observe that $\di R_Q\sim \di Q$.
Then multiply both sides with the same factor so that $(\di Q)^{-t}$
is the scale factor left for the LHS. Then 

$$
\eq{
||&\n^hu||_{H^{\l,pnt}(Q)}\cdot(\di Q)^{-t}\cdot
\G_{m,k,p,{\Cal A}}(x)^{1\over p}
\crd
&\le
A((\int_Q|\n^{k+1}u|^{p_1}dx)^{1\over p_1}
(\di Q)^{s_1\over p_1}
+(\int_{R_Q}|\n^mu|^pdx)^{1\over p}
(\di Q)^{s\over p}).
\cr}
\leqno(6.28)
$$

Since $\d^{-t}\sim (\di Q)^{-t}$ in the cube $Q$ we have simply

$$
\eq{
||\n^hu||&_{H^{\l,pnt}(Q,\L(x)^{1\over p}\d^{-t})}
\crd
\le&
A||\n^hu||_{H^{\l,pnt}(Q)}\cdot\G_{m,k,p,{\Cal A}}(x)^{1\over p}(\di Q)^{-t}.
\cr}
\leqno(6.29)
$$

Now raise (6.28) and (6.29) to power $\max\ \lbrace p,p_1\rbrace$.
Combine and sum over $Q$.
\ss

Then it holds

$$
||\n^hu||_{H^{\l,pnt}(\O,\L(x)^{1\over p}\d^{-t})}^{[p,p_1]}
\le
\sum_Q||\chi_Q\n^hu||_{H^{\l,pnt}(\O,\L(x)^{1\over p}\d^{-t})}
^{[p,p_1]}.
\leqno(6.30)
$$

This follows since the pointwise H\"older qoutient has local
properties and the contribution in the LHS of (6.30) is (roughly)
comming from one cube only and this contribution is found in (roughly)
one of cubes of the RHS of (6.30).
\ss

Next we we calculate the new RHS after these operations have been done.
\ms

Though first some (basic) inequalities are pointed out.
\ss

Note that by equivalent norms and/or quasinorms in
finite dimensions we have for $a,b,r>0$ that

$$
(a^r+b^r)^{1\over r}\le A(r)(a+b)
\leqno(6.31)
$$

There are also two well-known elementary inequalities for nonnegative
numbers.

\proclaim{Lemma 6.18}
$$
\text{For } r\ge 1\ \ \sum_{n=1}^{\infty}|a_n|^r
\le(\sum_{n=1}^{\infty} |a_n|)^r;
\ \ r\le 1\ \ \sum_{n=1}^{\infty} |a_n|^r \ge (\sum_{n=1}^{\infty} |a_n|)^r.
\leqno(6.32\,a,b)
$$
\endproclaim

Proof. Exercise.
\ms

We begin with (6.28), follow up the calculations and make use of
Lemma 6.18. Then

$$
\sum_Q(\int_Q|\n^{k+1}u|^{p_1}\d^{s_1}dx)^{[p,p_1]\over p_1}
\le
(\int_{\O}|\n^{k+1}u|^{p_1}\d^{s_1}dx)^{[p,p_1]\over p_1}.
\leqno(6.33)
$$

In the same way the second term of the same RHS  becomes

$$
\eq{
\sum_Q&(\int_{R_Q}|\n^mu|^pdx)^{[p,p_1]\over p}
(\di Q)^{{s\over p}{[p,p_1]}}
\crd
&\le
(\sum_Q\int_{R_Q}|\n^mu|^p\d^s\d^{-s}dx\cdot(\di Q)^s)^{[p,p_1]\over p}.
\cr}
\leqno(6.34)
$$

The last expression is estimated from above with Lemma 6.4. Then it holds

$$
\eq{
\sum_Q\int_{R_Q}|\n^mu|^p\d^s\d^{-s}dx
\cdot
(\di Q)^s
\crd
\le 
{A\over 1-2^{-s}}
\int_{\O}|\n^mu|^p\d^sdx.
\cr}
\leqno(6.35)
$$

The results are collected as

$$
\eq{
||\n^hu||&_{H^{\l,pnt}(\O,\L(x)7^{1\over p}\d^{-t})}^{[p,p_1]}
\le
A((\int_{\O}|\n^{k+1}u|^{p_1}\d^{s_1}dx)^{[p,p_1]\over p_1}
\crd
&+
({1\over 1-2^{-s}})^{[p,p_1]\over p})
(\int_{\O}|\n^mu|^p\d^sdx)^{[p,p_1]\over p}).
\cr}
\leqno(6.36)
$$

But by (6.31) we can take a root of these terms and finally obtain

$$
\eq{
||\n^hu||_{H^{\l,pnt}(\O,\L(x)^{1\over p}\d^{-t})}
\le&
A((\int_{\O}|\n^{k+1}u|^{p_1}\d^{s_1}dx)^{1\over p_1}
\crd
&+
{1\over (1-2^{-s})^{1\over p}}
(\int_{\O}|\n^mu|^p\d^sdx)^{1\over p}).
\cr}
\leqno(6.37)
$$

End of proof Case A -- (6.23).

\proclaim{Comment}
\vskip 0.1truecm
Note that the constants for the two terms in the RHS 
differs in dependency on $s$.
\endproclaim
\ss

Proof Case A -- (6.24)
\ss

The starting point is again a Poincar\'e inequality, see [WAN5], but now
with another LHS. This makes no major difference. We conclude that

$$
(\int_{\tilde R_Q}| u|^qdx)
^{1\over q}
\le {A\over\G_{m,k,p,{\Cal A}}(x)^{1\over p}}
((\int_{\tilde Q}|\n^{k+1}u|^{p_1}dx)
^{1\over p_1}
+(\int_{\tilde R_Q}|\n^{m}u|^pdx)^{1\over p})).
\leqno(6.38)
$$

Then the RHS of the wanted $(6.24)$ is got in exactly the same way as in the 
earlier proof of (6.23), i.e. 
dilate and rearrange like in the beginning of that proof.
Then sum in the same way as in the proof of (6.23).

The result is the RHS of (6.24).

However this procedure gives another LHS

$$
\sum_Q(\int_Q|u|^q\cdot\G_{m,k,p}(x)^{q\over p}dx)^{[p,p_1]\over q}.
\leqno(6.39)
$$

If $q\ge [p,p_1]$, then the desired result follows from Lemma 6.18.

Hence the first part of the statement in Case A -- (6.24) is proved.
\ms

Let instead $0<q<[p,p_1]$. Then the procedure is much the same.
Something else is needed instead of Lemma 6.18 though. Here the tool is
the H\"older inequality for sums.
\ss

Since $s'$ is defined as ${[p,p_1]\over q}$ and $f$ is defined with
$||f||_{l^s}=1$, we observe that 

$$
\eq{
(\sum_Q& f(Q)^s)^{1\over s}
(\sum_Q(\int_Q|u|^q\cdot\G_{m,k,p}(x)^{q\over p}dx)^{s'})^{1\over s'}
\crd
&\ge
\sum_Q f(Q)\int_Q|u|^q\cdot\G_{m,k,p}(x)^{q\over p}dx
\crd
&=
\int_{\O}|u|^qf(x)\G_{m,k,p}(x)^{q\over p}dx.
\cr}
\leqno(6.40)
$$

and then raise both sides to power $s'$. We obtain

$$
\eq{
(\sum_Q& f(Q)^s)^{s'\over s}
\sum_Q(\int_Q|u|^q\G_{m,k,p}(x)^{q\over p}dx)^{s'}
\ge
(\int_{\O}|u|^qf(x)\G_{m,k,p}(x)^{q\over p}dx)^{s'}
\crd
&=
(\int_{\O}|u|^qf(x)\G_{m,k,p}(x)^{q\over p}dx)^{{1\over q}\cdot [p,p_1]}.
\cr}
\leqno(6.41)
$$

as before the result follows by taking a root of each
term and use (6.31).
\ss

End of proof of Case A -- (6.24).
\ss

End of proof of Case A.
\bs

Proof of Case B
\ss

The proof is almost the same as the one for Case A.
The difference is in the beginning of the argument.

We want to prove (6.23). We begin with (6.26) with exponent $p_0$ 
instead of $p$, 

$$
\eq{
||\n^h u||_{H^{\l,pnt}(\tilde Q)}
&\le 
{A\over\G_{m,k,p_0,{\Cal A}}(x))^{1\over p_0}}
((\int_{\tilde Q}|\n^{k+1}u|^{p_1}dx)^{1\over p_1}
\crd
&+(\int_{\tilde R_Q}|\n^mu|^{p_0}dx)^{1\over p_0})
\le
\cr}
\leqno(6.42)
$$

Use the H\"older inequality on the second term in the RHS and then
(6.42) continues

$$
\eq{
\le
&{A\over\G_{m,k,p_0}(x)^{1\over p_0}}
((\int_{\tilde Q}|\n^{k+1}u|^{p_1}dx)^{1\over p_1}
\crd
&+(\int_{\tilde R_Q}|\n^mu|^p\d^{s+a}dx)^{1\over p}
(\int_{\tilde R_Q}\d^{-{s+a\over {p-p_0\over p_0}}}dx
)^{p-p_0\over pp_0})
\cr}
\leqno(6.43)
$$

for $a>0$. 
\ss

Now the last integral has a uniform bound independent of $Q$ if there
is a positive $a$ with

$$
s+a<{p-p_0\over p_0}(N-\text{dim}_{loc}(\partial \O,\O)).
\leqno(6.44)
$$

This follows from the given condition on $s$ and the
definition of $\dim_{loc}$.

Then the proof is completed in the same way as for Case A.
\ss

End of proof of Case B -- (6.23)
\ss

The Case B -- (6.24) is similar.
\ss

End of proof of Case B.
\bs

Proof of Case C and D
\ss

The proofs proceed in the same way as for Case A and B above, but with the 
difference that the $\T$-capacity is used instead of the $\G$-capacity.
\ss

End of proof of Case C and D.
\bs

Proof of Case E
\ss

It follows from Case A that for $\b$ small and positive

$$
\sum_0^{m-1}(\int{|\n^k u|^p\de^{-\b}\over\de^{(m-k)p}}dx)^{1\over p}
\le
{A\over\b^{1\over p}}
(\int|\n^mu|^p\de^{-\b}dx)^{1\over p}.
\leqno(6.45)
$$

Make a change of the dependent variable as follows

$$
u'=u\cdot\de^{\b'-\b\over p}.
\leqno(6.46)
$$

-- This transformation is clearly a set-isomorphism of 
$C^{\infty}_0(\O)$ to itself.
\ss

This change of variable is evaluated for the following expression.
The triangle inequality has also been used.

$$
\eq{
&(\int{|\n^k u'|^p\de^{-\b'}\over\de^{(m-k)p}}dx)^{1\over p}
\crd
&\le
A(\int{|\n^k u|^p\de^{-\b}\over\de^{(m-k)p}}dx)^{1\over p}
+
A{|\b'-\b|\over p}\sum_{r=0}^{k-1}
(\int{|\n^r u|^p\de^{-\b}\over\de^{(k-r)p}}dx)^{1\over p}.
\cr}
\leqno(6.47)
$$

Then sum over $k$. We obtain

$$
\sum_{k=0}^{m-1}
(\int{|\n^ku'|^p\de^{-\b'}\over\de^{(m-k)p}}dx)^{1\over p}
\le
A\sum_{k=0}^{m-1}
(\int{|\n^ku|^p\de^{-\b}\over\de^{(m-k)p}}dx)^{1\over p}.
\leqno(6.48)
$$

Next repeat (6.47) but with $m$ instead and use of the triangle inequality

$$
\eq{
(\int{|\n^mu'|^p\de^{-\b'}\over\de^{(m-k)p}}dx)^{1\over p}
\ge
&A'(\int{|\n^mu|^p\de^{-\b}\over\de^{(m-k)p}}dx)^{1\over p}
\crd
&-A''{|\b'-\b|\over p}
\sum_{k=0}^{m-1}
(\int{|\n^ku|^p\de^{-\b}\over\de^{(m-k)p}}dx)^{1\over p}.
\cr}
\leqno(6.49)
$$

However if say

$$
{A'\over \b^{1\over p}}\ge 2{A''|\b'-\b|\over p},
\leqno(6.50)
$$

then the first term of the RHS of (6.49) dominates and the
inequality (6.24) with only the higher order term in RHS holds with
$s=\b'$ say.

It remains to calculate the possible $\b'$:s. We make the choice $\b'=-\b$ 
in order to simplify. Then $\b$ has to be chosen so that

$$
c
\ge
\b^{1-{1\over p}},
\leqno(6.51)
$$

for some constant $c>0$, which can be calculated.
This is possible to do if $p>1$ and $\b$ small enough.
\ss

End of proof of Case E.
\ss

End of proof of Theorem 6.17.
\bs

We will give examples in Part II of how the two-integral RHS in
Theorem 6.17 in certain cases can be used to get a better 
one-integral in the RHS Hardy inquality
than the one-integral formulation taken directly from 
Theorem 6.17.

-- This is the motivation for the two-integral RHS
formulation in Theorem 6.17.
\ss

Next we give a corollary to Theorem 6.17. 

It consists of a list of cases where
the one-integral RHS of 6.17 have been used.
The weights are specified to give a dilation invariant formulation.
Hence the weights are similar to the original one dimensional Hardy inequality.

We include the $f$-cases with $q<p$ since they can be stated simultanuously.
\ss

\proclaim{Corollary 6.19 to Theorem 6.17} For notation see 6.16.
\ss

Let $u\in {\Cal A}$
\ss

In the cases (i)-(x) it follows from [WAN5] that there is a $b>0$ with

$$
\G_{m,m-1,{\Cal A}}(x)\ge b>0.
\leqno(6.52)
$$
\ss

Preconditions for (6.55).

$$
\cases
(m-h)p>N>(m-h-1)p,
\\
0<l\le m-h-{N\over p},
\\
t=m-h-l-{s\over p}-{N\over p}.
\endcases
\leqno(6.53)
$$
\ss

Preconditions for (6.56).

$$
\cases
t=mq-({q\over p}-1)N-{q\over p}s
\\
and
\\
N>mp$ with $0<q\le{pN\over N-mp}
\\
or
\\
N\le mp$ with $0<q\le\infty.
\endcases
\leqno(6.54)
$$

$$
||\n^hu||_{H^{\l,pnt}(\O,\d^{-t})}
\le
A(\int|\n^mu|^p\d^sdx)^{1\over p}.
\leqno(6.55)
$$

$$
(\int{| u|^qF(x)\over\d^t}dx)^{1\over q}
\le
A(\int|\n^mu|^p\d^sdx)^{1\over p}.
\leqno(6.56)
$$

Here put $F(x)=1$ if $q\ge p$ and $F(x)=f(x)$ if $q<p$.
\ms

Then (6.55) respectively (6.56) holds in cases (i)-(x).
\ss

To be added:

For {\rm(ii)} and {\rm(vi)} let $1\le p_0<p$,

for {\rm(iv)}, {\rm(viii)} and {\rm(x)} let $1<p_0<p$,

for {\rm (v)-(viii)}, let there exist an $r$, $0\le r\le N$, 
such that for all $Q\in {\Cal F}_{\O}$ there is set of
orthogonal projections onto hyperplanes, $\{S_i\}^{N-r}_{i=1}$, 
with $r$-dim cube $Q'$

$$
Q'\subset (\prod_{i=1}^{N-r}S_i)\ (\tilde R_Q\cap\tilde\O^c)
$$ 

and $\di Q'\ge b>0$ for all $Q\in{\Cal F}_{\O}$.
\ms

Below $Q\in{\Cal F}_{\O}$ and $A$, $s_0$ are positive constants
independent of $u$.
\ms

\roster
\item"(i)" 
Let ${\Cal A}=W^{m,p}_0(\O)$ and let 
$C_{1,p}(\tilde \O^c\cap \tilde R_Q)\ge \text{\rm const.}>0$,

If $s<0$, then $p\ge 1$ and if $0\le s<s_0$, then $p>1$.
\endroster
\ss

\roster
\item"(ii)" 
Let ${\Cal A}=W^{m,p}_0(\O)$ and let 
$C_{1,{p_0}}(\tilde \O^c\cap\tilde R_Q)\ge \text{\rm const.}>0$.
Let 

$s<({p\over p_0}-1)(N-\dim_{loc}(\pa\O,\O))$.
\endroster
\ss

\roster
\item"(iii)" 
Let ${\Cal A}=W^{2,p}_0(\O)_+$ for $p\ge 1$. 
Let 
$C_{2,p}(\tilde \O^c\cap \tilde R_Q)\ge \text{\rm const.}>0$.

If $s<0$, then $p\ge 1$, and if $0\le s<s_0$, then $p>1$.
\endroster
\ss

\roster
\item"(iv)" 
Let ${\Cal A}=W^{2,p}_0(\O)_+$ and $p>1$. 
Let $C_{2,{p_0}}(\tilde \O^c\cap\tilde R_Q)\ge \text{\rm const.}>0$.
Let $s<({p\over p_0}-1)(N-\dim_{loc}(\pa\O,\O))$.
\endroster
\ss

\roster
\item"(v)" 
Let ${\Cal A}=W^{m,p}_0(\O)$ and let $p>N-r$. Let $s<s_0$ and let $p>1$.
\endroster
\ss

\roster
\item"(vi)" 
Let ${\Cal A}=W^{m,p}_0(\O)$, let $p_0>N-r$
and let $s<({p\over p_0}-1)(N-\dim_{loc}(\pa\O,\O))$.
\endroster
\ss

\roster
\item"(vii)" 
Let ${\Cal A}=W^{2,p}_0(\O)_+$ and let $2p>N-r$.
Let $s<s_0$ and let $p>1$.
\endroster
\ss

\roster
\item"(viii)" 
Let ${\Cal A}=W^{2,p}_0(\O)_+$, let $2p_0>N-r$,
and let 
\pn
$s<({p\over p_0}-1)(N-\dim_{loc}(\pa\O,\O))$.
\endroster
\ss

\roster
\item"(ix)" 
Let ${\Cal A}=W^{m,p}_0(\O)$, let $s<0$, let $p>1$,
let $C_{m,p}(\O^c)\not=0$,
let $\O^c$ be
selfsimilar and let $\pa\O$ not be a subset of a hyperplane.
\endroster
\ss

\roster
\item"(x)" 
Let ${\Cal A}=W^{m,p}_0(\O)$,
let $s<({p\over p_0}-1)(N-\dim_{loc}(\pa\O,\O))$,
let 
\pn
$C_{m,p_0}(\O^c)\not=0$,
let $\O^c$ be selfsimilar 
and let $\pa\O$ not be a subset of a hyperplane.
\endroster
\endproclaim
\ms

{\bf Some earlier results.}
\ss

Ancona has treated the case $q=p$, $p>N$ and $\O$ bounded in [ANC1]
and the case (i) in Corollary 6.19 for $m=1$, $p=2$, $s=0$ in [ANC2].
Lewis has treated the case (i) with $m=1$, general $p$, $s\le 0$ in [LEW].
\ss

In [WAN2] it is given a deliberately short version, with a different
proof. (Here is also reference made to more extensive unpublished 
work material).

The results of Ancona in [ANC1], [ANC2] and Lewis in 
[LEW] on Hardy inequalities are special cases of the result in [WAN2] ,
except their necessity condition.

This have been discussed in the introductory part.
\ms

{\bf Observation.}
\ss

In Theorem 6.17 the constant $A$ does not depend on $q$ in (6.25).
This property has been the subject of a paper by Kavian, see [KAV],
whose result was used by Brezis and Turner
on a nonlinear partial differential equation and the corresponding
eigenvalue problem, see [BRE--TUR].
\ms

The paper [BRE-BRO] by Brezis and Browder got
the following problem into focus. The question is, when holds

$$
W^{m,p}_0(\O)=W^{m,p}_0(\O)_+-W^{m,p}_0(\O)_+.
\leqno(6.58)
$$

Ancona gave some answer ($p>1$), see [ANC1], and
he proved that, if either $\O$ was bounded Lipschitz or if $\O$ was bounded 
and $p>N$, then (6.58) holds.

\bs

The following result is given here

\proclaim{Theorem 6.20}
Let one of the conditions {\rm(i)}, {\rm(ii)}, {\rm(v)},
{\rm(vi)}, {\rm(ix)} and {\rm(x)} of Corollary 6.19 hold together with
eventual extra conditions given in the Corollary.
Let $q=p$.
Then 

$$
W^{m,p}_0(\O,\d^s)=W^{m,p}_0(\O,\d^s)_+-W^{m,p}_0(\O,\d^s)_+,
\leqno(6.59)
$$
also the $W^{m,p}(\O,\d^s)$-norms 
for the nonnegative functions are less than a factor times the 
norm of the original function. ( -- A feature also in the Ancona result.)
\endproclaim

Proof.
\ss

Follows from Theorem 6.17 and Theorem 6.23.
\bs

There is another interesting result.
\ss

\proclaim{Theorem 6.21}
Let $u\in W^{2,p}_0(\O,\d^s)$ and $1<p$,
If one of the conditions (i)--(iv) below is satisfied, then 

$$
u\in W^{2,p}_0(\O,\d^s)_+-W^{2,p}_0(\O,\d^s)_+
\leqno(6.60)
$$

holds, if and only if

$$
\int|u|^p\d^{-2p+s}dx<\infty.
\leqno(6.61)
$$
\ss

Additional comments:
\ss

In {\rm (iv)} 
let there exist an $r$, $0\le r\le N$, such that to all
$Q\in {\Cal F}_{\O}$ there is set of
orthogonal projections onto hyperplanes, $\{S_i\}^{N-r}_{i=1}$, 
with $r$-dim cube $Q'$

$$
Q'\subset (\prod_{i=1}^{N-r}S_i)\ (\tilde R_Q\cap\tilde\O^c)
$$ 

and $\di Q'\ge b>0$ for all $Q\in{\Cal F}_{\O}$.
\ss

Below $Q\in{\Cal F}_{\O}$ and $A$, $s_0$ are positive constants
independent of $u$.)
\ss

\roster
\item"(i)" 
Let 
$C_{2,p}(\tilde \O^c\cap \tilde R_Q)\ge \text{\rm const}>0$ for all
$Q\in {\Cal F}_{\O}$.
If $s<0$, then $p\ge 1$,
and if $s<s_0$, then $p>1$.
\endroster
\ss

\roster
\item"(ii)" 
Let 
$C_{2,{p_0}}(\tilde \O^c\cap\tilde R_Q)\ge \text{\rm const}>0$ for all
$Q\in {\Cal F}_{\O}$.
Let $s<({p\over p_0}-1)(N-\dim_{loc}(\pa\O,\O))$.
Let $1\le p_0<p$.
\endroster
\ss

\roster
\item"(iii)" 
Let $2p>N-r$, let $s<s_0$ and let $p>1$.
\endroster
\ss

\roster
\item"(iv)" 
Let $2p_0>N-r$ and let 
$s<({p\over p_0}-1)(N-\dim_{loc}(\pa\O,\O))$. Let $1<p_0<p$.
\endroster
\endproclaim
\bs

Proof.
\bs

The sufficiency follows from Theorem 6.23. The neccessity follows from 
Corollary 6.19.
\bs

The following theorem can be read out from the proof 
of the first theorem by Ancona in [ANC1]. 
We give a different theorem though.
Since Ancona's account is brief we repeat the arguments andz more generously.

-- This version is more general.

\proclaim{Theorem 6.22}
Let $1<p$ and let $s$ be a real number. Let $\O$ be open subset of \R,
with non-empty complement.
Let
$$
u\in W^{m,p}_0(\O,\d^s).
$$

Then 
$$
\int_{\O}|u|^p\d^{-mp+s}dx<\infty
$$
implies that there exist
$$
u_i\in W^{m,p}_0(\O,\d^s)_+,\ i=1,2,
$$
with $u=u_1-u_2$.
\endproclaim

Proof.
\ss

Let $u_Q\in W_0^{m,p}(\a Q)$, with say $\a=4/3$.

If there exists $v_Q\in W_0^{m,p}(\b Q)$, with say $\b=4^2/3^2$,
such that $v_Q\ge u_Q$ and $v_Q\ge 0$ and if furthermore it
holds that 

$$
||v_Q||_{W^{m,p}}\le A_0||u_Q||_{W^{m,p}},
$$
then the problem is solved locally in a sense.

The procedure now is first to show that this implies that there is
a global solution. After that we return and solve the local problem.
\ss

Let $Q_0$ be a unit cube. 
Take an $\eta\in C^{\infty}_0({4\over 3}Q_0)_+$, $\eta|_{Q_0}=1$. 
Let $Q$ be any cube, centre $x_Q$. Define

$$
\eta_Q=\eta({x-x_Q\over l(Q)})
\leqno(6.62)
$$
Let $Q\in {\Cal F}_{\O}$, a set of Whitney cubes for $\O$.
Let $u\in W_0^{m,p}(\O,\d^s)$. Define $u_Q=\eta_Qu$, then by
assumption there is an $v_Q$ with properties as above.
Define

$$
v=\sum_{Q\in{\Cal F}_{\O}}v_Q.
\leqno(6.63)
$$
Clearly $v\ge 0$ and $v\ge u$.

In this global situation it only remains to make the norm estimate,

$$
\eq{
||v||^p_{W^{m,p}(\O,\d^s)}&\le A\sum_{k=0}^m||\n^kv||^p_{L^p(\O,\d^s)}
\crd
&=
A\sum_{k=0}^m
||\sum_{Q\in{\Cal F}_{\O}}\n^kv_Q||^p_{L^p(\O,\d^s)}
\le
\cr}
\leqno(6.64)
$$

Next we use the fact that the enlarged cubes only have at most
a fixed number of overlap. This means that in this situation
the integrals can be decomposed at the cost of a constant factor only,

$$
\eq{
&\le
A\sum_{k=0}^m
\sum_{Q\in{\Cal F}_{\O}}
||\n^kv_Q||^p_{L^p(\O,\d^s)}
\crd
&\le A\sum_{k=0}^m
\sum_{Q\in{\Cal F}_{\O}}
||\n^kv_Q||^p_{L^p(\O)}l(Q)^s
\crd
&\le 
A\sum_{Q\in{\Cal F}_{\O}}
||v_Q||^p_{W^{m,p}(\O)}l(Q)^s
\crd
&\le
A\sum_{Q\in{\Cal F}_{\O}}
||u_Q||^p_{W^{m,p}(\O)}l(Q)^s
\le
\cr}
\leqno(6.65)
$$
Next step is a standard interpolation of norm and the scale factor
involved is e.g. the result of a dilation of a unit dimension case

$$
\eq{
&\le
A\sum_{Q\in{\Cal F}_{\O}}
||u_Q||^p_{L^p(\O)}l(Q)^{-mp+s}
+
||\n^mu_Q||^p_{L^p(\O)}l(Q)^s
\crd
&\le
||u||^p_{L^p(\O,\d^{-mp+s})}
+
||\n^mu||^p_{L^p(\O,\d^s)}.
\cr}
\leqno(6.66)
$$

The last step is a consequence of that the $p$-powers make norms into
integrals and then this inverse kind of inequality becomes natural.
\ss

In (6.66) it is clear that the first term in the last expression is finite
and that then $u\in W^{m,p}(\O,\d^s)$ makes the
whole expression finite. Hence the wanted function $v$ exist with
finite norm. This step-wise build-up of $v$ also ensures that
$v\in W_0^{m,p}(\O,\d^s)$ as wanted.

Hence it remains only to prove the local cube-wise property.

Let $u\in W^{m,p}(\O,\d^s)$ and $u_Q=\eta_Q u$. 
The Bessel kernel $G_m$ together with an $f\in L^p$ solves

$$
G_m*f=u_Q.
\leqno(6.67)
$$

This is $u_Q$ written as a Bessel potential, ($1<p$),
and the norm $||f||_{L^p}$ is equivalent to the Sobolev space norm.
Further $G_m$ acts positively, i.e. if $0\le g$ and $g\in L^p$ then
$0\le G_m*g$. In the situation above take $f_+$ instead of $f$ and
then $G_m*f\le G_m*f_+$ and of course $||f_+||_{L^p}\le ||f||_{L^p}$.

Note that as usual the Sobolev space $W^{m,p}$ is defined
by the quasicontinuous functions only. 

To ensure the short range of support of $v_Q$ we define 

$$
v_Q=\eta_{(4^2/3^2Q)}(G_m*f_+).
\leqno(6.68)
$$

By the standard Poincar\'e inequality and the construction of the Bessel
potentials it holds that

$$
||v_Q||_{W^{m,p}}\le ||f||_{L^p}=||G_m*f||_{W^{m,p}}=||u_Q||_{W^{m,p}}.
\leqno(6.69)
$$

Hence the local problem is solved.
\ss

The proof is complete.

\enddocument